\documentclass[11pt]{amsart}
\usepackage[T1]{fontenc}
\usepackage[utf8]{inputenc}
\usepackage[english]{babel}
\usepackage{amsmath,amssymb,amsthm,mathtools,mathrsfs}
\usepackage[a4paper,left=28mm,right=28mm,top=25mm,bottom=25mm]{geometry}
\usepackage{lmodern}
\usepackage{microtype}
\usepackage[colorlinks=true,allcolors=blue]{hyperref}
\hypersetup{pdftitle={Diophantine definitions of nonsplitting},pdfauthor={Guang Hu}}

\newtheorem{theorem}{Theorem}[section]
\newtheorem{proposition}[theorem]{Proposition}
\newtheorem{lemma}[theorem]{Lemma}
\newtheorem{corollary}[theorem]{Corollary}
\theoremstyle{definition}
\newtheorem{definition}[theorem]{Definition}
\newtheorem{remark}[theorem]{Remark}
\DeclareMathOperator{\Br}{Br}
\DeclareMathOperator{\Trd}{Trd}
\DeclareMathOperator{\Nrd}{Nrd}
\DeclareMathOperator{\ind}{ind}
\DeclareMathOperator{\per}{per}
\DeclareMathOperator{\inv}{inv}
\DeclareMathOperator{\SB}{SB}

\DeclareMathOperator{\End}{End}
\DeclareMathOperator{\chara}{char}
\newcommand{\CO}{\mathcal O}

\newcommand{\Km}{\mathfrak m}
\newcommand{\Cc}{\mathcal C}
\newcommand{\Nc}{\mathcal N}
\newcommand{\op}{\mathrm{op}}
\newcommand{\A}{\mathbb A}
\newcommand{\Q}{\mathbb Q}
\newcommand{\Z}{\mathbb Z}
\newcommand{\F}{\mathbb F}
\newcommand{\GL}{\operatorname{GL}}
\newcommand{\Norm}{\operatorname{N}}

\title[Diophantine definitions of nonsplitting]
  {Diophantine definitions of nonsplitting for central simple algebras
   over global fields}
\author{Guang Hu}
\address{The Mathematics School, Shandong University, 27 Shanda Nanlu,
Jinan, Shandong, P. R. China}
\address{Faculty of Mathematics, Technische Universit\"at Dresden,
01062 Dresden, Germany}
\email{202490000031@sdu.edu.cn}
\email{huguang@mail.ustc.edu.cn}
\email{guang.hu@mailbox.tu-dresden.de}
\date{}

\begin{document}
\begin{abstract}
For every global field and every fixed degree, we prove that the
structure constants defining nonsplit central simple algebras form a
diophantine set. We formulate a local witness principle separating
reliable localization, exact localization, and local detection.
Adapted Kummer and Artin--Schreier presentations provide a common
construction of witnesses through a uniform radical formula and
ramification matching. This gives an independent proof
of the non-characteristic case, which also follows from Morrison's
nonnorm theorem, and treats the characteristic-degree case by the same
method. We recover prime-degree nonnorm definability and give applications
to rootless polynomials, anisotropic quadratic forms in all
characteristics, and generalized Severi--Brauer varieties.
\end{abstract}
\maketitle

\section{Introduction}
Let $K$ be a global field. A subset of $K^n$ is called
\emph{diophantine over $K$} if it is defined by a finite system of
polynomial equations with coefficients in $K$ and existentially
quantified variables. We allow parameters from $K$. In other words, we have the following definition.
\begin{definition}\label{def:diophantine}
Let $X\subseteq K^n$. We say that $X$ is \emph{diophantine over $K$}
if there exist an integer $r\geq 0$ and a polynomial
\[
f\in K[x_1,\ldots,x_n,y_1,\ldots,y_r]
\]
such that, for every $\mathbf{x}\in K^n$,
\[
\mathbf{x}\in X
\quad\Longleftrightarrow\quad
\exists \mathbf{y}\in K^r
\quad f(\mathbf{x},\mathbf{y})=0.
\]
\end{definition}
For a family of subsets parameterized by elements of $K^m$, we define
uniform diophantineness as follows.

\begin{definition}\label{def:uniformly-diophantine}
Let $P\subseteq K^m$, and let
\[
\{X_{\mathbf t}\}_{\mathbf t\in P},
\qquad X_{\mathbf t}\subseteq K^n,
\]
be a family of subsets. We say that the family is
\emph{uniformly diophantine over $K$} if the relation
\[
\{(\mathbf t,\mathbf x)\in K^{m+n}:
\mathbf t\in P,\ \mathbf x\in X_{\mathbf t}\}
\]
is diophantine over $K$.
\end{definition}
Finite unions, finite intersections, and projections of diophantine sets
are diophantine. There is no general closure theorem under complements.
Nevertheless, several natural arithmetic conditions and their negations
both admit diophantine definitions. Examples include reducibility and
irreducibility of polynomials \cite{Dittmann2018}, isotropy and anisotropy
of quadratic forms in characteristic different from two \cite{hu2026},
and membership and non-membership in certain norm groups
\cite{KTU,Morrison2019}. We study the corresponding question for
central simple algebras. Splitting is diophantine, since it is expressed
by an isomorphism with a matrix algebra; our main result concerns
nonsplitting.

Fix $n\ge2$ and put $N=n^2$. A multiplication table on $K^N$ is
specified by
\[
 \mathbf c=(c_{ij}^k)_{1\le i,j,k\le N},
 \qquad e_i e_j=\sum_{k=1}^N c_{ij}^k e_k.
\]
Let $\Cc_n(K)$ be the set of tables defining central simple
$K$-algebras with identity $e_1$. For $\mathbf c\in\Cc_n(K)$, write
$A_{\mathbf c}$ for the resulting algebra. The set $\Cc_n(K)$ is
itself diophantine (see Lemma~\ref{csa locus}).
Our main theorem is the following:
\begin{theorem}\label{main theorem}
Let $K$ be a global field and $n\ge2$. Then
\[
 \Nc_n(K)=\{\mathbf c\in\Cc_n(K):
                  A_{\mathbf c}\text{ is nonsplit over }K\}
\]
is diophantine over $K$.
\end{theorem}

Period--index equality over global fields reduces this statement to
prime degree (see Proposition~\ref{prop:prime-reduction}):
\begin{theorem}\label{prime version}
Let $K$ be a global field and let $\ell$ be a prime. Then $\Nc_\ell(K)$
is diophantine over $K$.
\end{theorem}

The prime-degree case away from the characteristic can already be
deduced from Morrison's theorem on cyclic norm forms
\cite[Theorem~1.5]{Morrison2019}, together with the standard theory of
central simple algebras and the diophantine encoding of cyclic
extensions. We give an independent proof of this case. Our purpose is
to organize the local witness mechanism appearing in earlier work into
a common construction, and to apply the same construction in the
characteristic-degree case, which is outside the range of that theorem.

The proof separates three tasks: ensuring that an auxiliary parameter
has at least one genuine witness place, choosing a parameter that
isolates any prescribed place, and detecting the desired local
obstruction. The general principle is stated in
Proposition~\ref{prop:general-witness}; its Brauer-theoretic realization
is Proposition~\ref{local-diophantine-correspondence}. The latter uses
the cyclicity of the local prime-torsion Brauer group to replace a
nonvanishing condition by a splitting certificate for an auxiliary
algebra.

The main construction uses adapted Kummer presentations away from the
characteristic and adapted Artin--Schreier presentations in the
characteristic degree. Both yield the same valuation condition, so a
single radical formula and a single ramification-matching argument apply to
both. Ramification matching fixes the difference of two auxiliary
Brauer classes, forcing their ramification sets to agree outside two
prescribed places; reciprocity ensures that their common ramification
set is nonempty. Places occur in the proof but are never quantified in the
defining formula.

The semilocal rings used here are due to Dittmann
\cite[Proposition~2.9]{Dittmann2018}. The use of parameterized semilocal
rings and their radicals to encode local obstructions is already
present in \cite{KTU,Morrison2019}. Related uniform constructions,
including characteristic two, were developed by Daans
\cite{Daans2021} for universal definitions of subrings.
Our contribution here is the organization
of these tasks into the witness principle and the common construction
above. In particular, our proof uses neither Morrison's nonnorm
definability theorem nor Dittmann's irreducibility theorem. A comparison
with earlier proofs, including applications to nonintegrality and
rootlessness, is given in Section~\ref{sec:applications}.

For a central simple algebra $A$ of degree $n$ and $1\le r<n$, let
$\SB_r(A)$ denote the variety of right ideals of reduced dimension $r$.
In particular, $\SB_1(A)$ is the usual Severi--Brauer variety.
We obtain the following application:
\begin{corollary}\label{cor:specialization}
Let $K$ be a global field, let $A$ be a central simple algebra of degree
$n$ over $K(T)$, and let $1\le r<n$. There is a nonempty open subset
$U\subseteq\A^1_K$ on which $A$ specializes to central simple
$K$-algebras $A_t$ of degree $n$, and the following set is diophantine:
\[
 \{t\in U(K):\SB_r(A_t)(K)=\varnothing\}.
\]
\end{corollary}

The proof uses the equivalence
$\SB_r(A)(K)\ne\varnothing\Longleftrightarrow\ind(A)\mid r$.
Over a global field, $\ind(A)=\per(A)$, so the failure of this
condition is exactly the nonsplitting of $A^{\otimes r}$.

\section{Algebraic and diophantine preliminaries}

\subsection{Central simple algebras over global fields}
We use the standard theory of central simple algebras, cyclic algebras,
and Brauer groups; see \cite[Chapters~2, 4, and~6]{GilleSzamuely2017}.
\begin{definition}
    The \emph{index} of $A$, denoted by $\operatorname{ind}(A)$,
is the degree of the division algebra $D$ in the Brauer class of $A$;
see \cite[Definition~4.5.1]{GilleSzamuely2017}.
The \emph{period} of $A$, denoted by $\operatorname{per}(A)$,
is the order of $[A]$ in $\operatorname{Br}(K)$;
see \cite[Definition~4.5.12]{GilleSzamuely2017}.
\end{definition}

\begin{lemma}
\label{lem:period-index}
Let $A$ be a central simple algebra over a field $K$.

\begin{enumerate}
    \item One has
    \[
    \operatorname{per}(A)\mid\operatorname{ind}(A)\mid\deg(A).
    \]

    \item If $K$ is a global field, then
    \[
    \operatorname{per}(A)=\operatorname{ind}(A),
    \]
    and every central division algebra is cyclic.
\end{enumerate}
\end{lemma}

\begin{proof}
The divisibility $\per(A)\mid\ind(A)$ is
\cite[Proposition~4.5.13]{GilleSzamuely2017}. Writing
$A\simeq M_r(D)$ with $D$ division gives
$\deg(A)=r\deg(D)=r\ind(A)$, proving the other divisibility.
For a global field, period--index equality and cyclicity follow from
the global Brauer theory recalled in
\cite[Remarks~6.5.5 and~6.5.6]{GilleSzamuely2017}.
\end{proof}

Let $V_K$ be the set of nonarchimedean places of a global field $K$,
and let $V_K^{\mathrm r}$ be its set of real places. We use the
contracted valuation ring and maximal ideal
\[
 \CO_v=\{x\in K:v(x)\ge0\},\qquad
 \Km_v=\{x\in K:v(x)>0\},
\]
where $v:K^\times\to\Z$ is normalized and $v(0)=+\infty$.
Write
\[
 \Delta(A)=\{v\in V_K:A\otimes_K K_v\text{ is nonsplit}\}.
\]
Thus $\Delta(A)$ records only finite places. We use Brauer reciprocity
in the form
\begin{equation}\label{eq:reciprocity}
 0\longrightarrow\Br(K)\longrightarrow
 \bigoplus_{v\text{ a place of }K}\Br(K_v)
 \xrightarrow{\ \sum_v\inv_v\ }\Q/\Z\longrightarrow0.
\end{equation}
At a finite place, $\Br(K_v)\simeq\Q/\Z$; at a real place its image
under the invariant map is $\{0,1/2\}$, and at a complex place the
Brauer group is zero. In particular, invariants $1/\ell,-1/\ell$ at
two distinct finite places and zero elsewhere give a degree-$\ell$
algebra, split at all real places. Local invariants are multiplied by
the local degree under finite extension. See
\cite[Section~6.5]{GilleSzamuely2017}.

\subsection{Encoding algebra operations}
\begin{lemma}\label{lem:algebra-isomorphism}
For fixed dimensions, isomorphism of unital algebras given by structure
constants is a diophantine relation. On central simple algebras of
fixed degrees, equality of Brauer classes is also diophantine.
\end{lemma}
\begin{proof}
Suppose $e_i e_j=\sum_k a_{ij}^k e_k$ and
$f_i f_j=\sum_k b_{ij}^k f_k$. A linear map with matrix
$X=(x_{ri})$, taking $e_i$ to $\sum_r x_{ri}f_r$, preserves
multiplication exactly when
\[
 \sum_k a_{ij}^k x_{sk}
   =\sum_{r,t}x_{ri}x_{tj}b_{rt}^s
 \quad\text{for all }i,j,s.
\]
Preservation of the specified identity is polynomial, and invertibility
is expressed by $z\det X=1$ for an additional variable $z$.

If $A$ and $B$ have degrees $a$ and $b$, respectively, then
\begin{equation}\label{eq:brauer-equality}
 [A]=[B]\quad\Longleftrightarrow\quad
 M_b(A)\simeq M_a(B).
\end{equation}
Indeed, the two matrix algebras have the same degree, and central simple
algebras of equal degree in the same Brauer class are isomorphic.
Tensor products, opposite algebras, and matrix algebras have structure
constants polynomial in the original tables. This proves the claim.
For matrix algebras, one may use a fixed change of basis containing the
identity if a table in the identity-first convention is required.
\end{proof}

\begin{lemma}\label{csa locus}
The set $\Cc_n(K)$ is diophantine over every field $K$.
\end{lemma}
\begin{proof}
The identity and associativity conditions are
\[
 c_{1i}^k=c_{i1}^k=\delta_{ik},\qquad
 \sum_q c_{ij}^q c_{qk}^s=\sum_q c_{jk}^q c_{iq}^s.
\]
For an associative unital algebra $A$, central simplicity is equivalent
to invertibility of the canonical map
\[
 A\otimes_K A^{\op}\longrightarrow\End_K(A),
 \qquad a\otimes b\longmapsto(x\mapsto axb).
\]
Its matrix has entries
$R_{(s,r),(i,j)}=\sum_q c_{ir}^q c_{qj}^s$.
Thus the remaining condition is $z\det R=1$ for some $z\in K$.
\end{proof}
Let $A$ be a central simple $K$-algebra of degree $n$.
For $x\in A$, write its reduced characteristic polynomial as
\[
P_x(Z)=Z^n+d_1Z^{n-1}+\cdots+d_n\in K[Z].
\]
This is the characteristic polynomial of the matrix representing
$x$ under an isomorphism $A\otimes_K L\simeq M_n(L)$ over a
splitting field $L$; it is independent of these choices.
The reduced trace and reduced norm of $x$ are
\[
\Trd_A(x)=-d_1,
\qquad
\Nrd_A(x)=(-1)^n d_n,
\]
respectively.
\begin{lemma}\label{lem:reduced-functions}
Fix $n\geq 2$. The set
\[
\left\{
(\mathbf c,\mathbf x,t,u)
\in \mathcal C_n(K)\times K^{n^2}\times K^2 :
t=\Trd_{A_{\mathbf c}}(x),\
u=\Nrd_{A_{\mathbf c}}(x)
\right\},
\qquad
x=\sum_{i=1}^{n^2}x_i e_i,
\]
is diophantine over $K$.
\end{lemma}
\begin{proof}
For $\mathbf c\in\mathcal C_n(K)$, let $A=A_{\mathbf c}$ and
$x=\sum_{i=1}^{n^2}x_i e_i$. The matrix $M_{\mathbf c,\mathbf x}$
of left multiplication by $x$ in the basis
$e_1,\ldots,e_{n^2}$ has entries
\[
(M_{\mathbf c,\mathbf x})_{kj}
=\sum_{i=1}^{n^2}x_i c_{ij}^k.
\]
In particular, its entries are polynomial in $\mathbf c$ and
$\mathbf x$.

Let
\[
P_x(Z)=Z^n+d_1Z^{n-1}+\cdots+d_n
\]
be the reduced characteristic polynomial of $x$. Then
\[
\det(ZI_{n^2}-M_{\mathbf c,\mathbf x})=P_x(Z)^n.
\]
Indeed, after extension to a splitting field, $A$ becomes $M_n$,
and left multiplication by $x$ acts on each of the $n$ columns
by the matrix representing $x$.

Moreover, $P_x$ is the unique monic polynomial of degree $n$
satisfying this identity. To see this, if two monic polynomials
$P,Q\in K[Z]$ satisfy $P^n=Q^n$, unique factorization in $K[Z]$
implies $P=Q$. This argument is valid even when
$\operatorname{char}(K)$ divides $n$.

Consequently, for $\mathbf c\in\mathcal C_n(K)$, the equalities
\[
t=\Trd_A(x),\qquad u=\Nrd_A(x)
\]
hold if and only if there exist $d_1,\ldots,d_n\in K$ such that
\[
\det(ZI_{n^2}-M_{\mathbf c,\mathbf x})
=\bigl(Z^n+d_1Z^{n-1}+\cdots+d_n\bigr)^n,
\qquad
t=-d_1,\qquad u=(-1)^n d_n.
\]
Here the polynomial identity in $Z$ is expressed by equating
coefficients. For fixed $n$, this gives a finite system of
polynomial equations in
$\mathbf c,\mathbf x,t,u,d_1,\ldots,d_n$.
Since $\mathcal C_n(K)$ is diophantine by Lemma~\ref{csa locus},
existentially quantifying $d_1,\ldots,d_n$ proves the claim.
\end{proof}

\section{Semilocal rings and local splitting certificates}
\label{section4}\label{local-to-diophantine-section}
Fix a prime $\ell$. Throughout the next three sections, algebras
of degree $\ell$ are assumed to split at every real place. This is
automatic unless $K$ is a number field and $\ell=2$.
Let $\Cc_\ell^0(K)$ denote this class of tables. Its diophantine
definability, including the exceptional case, will be established in
Lemma~\ref{lem:real-places} without using nonsplitting definability.

For such an algebra $A$, put
\[
 S_A=\{\Trd_A(x):x\in A,\ \Nrd_A(x)=1\},\qquad
 T_A=\begin{cases}S_A,&\ell>2,\\ S_A-S_A,&\ell=2.\end{cases}
\]
Here $S_A-S_A$ denotes the set of pairwise differences.

\begin{proposition}[Dittmann]\label{dittmann semilocal ring}
One has
\begin{equation}\label{eq:semilocal-ring}
 T_A=\bigcap_{v\in\Delta(A)}\CO_v.
\end{equation}
Membership in $T_A$ is uniformly diophantine in the structure constants
of $A$. An empty intersection is understood to be $K$.
\end{proposition}
\begin{proof}
The equality is \cite[Proposition~2.9]{Dittmann2018}. The hypothesis at real places is exactly the one required there. Uniform
diophantine definability follows directly from the definition of $T_A$
and Lemma~\ref{lem:reduced-functions}.
\end{proof}

\begin{lemma}\label{lem:sum-intersections}
For finite subsets $E,F\subset V_K$, put
$\CO_E=\bigcap_{v\in E}\CO_v$ and
$J_E=\bigcap_{v\in E}\Km_v$, with empty intersections equal to $K$.
Then
\[
 \CO_E+\CO_F=\CO_{E\cap F},\qquad J_E+J_F=J_{E\cap F}.
\]
\end{lemma}
\begin{proof}
The forward inclusions are immediate. If $E$ or $F$ is empty,
both equalities follow from the convention on empty intersections.

For the reverse inclusion in the second equality, let
$x\in J_{E\cap F}$. By weak approximation, choose $y\in K$
such that
\[
 y\in\Km_v \quad (v\in E),\qquad
 x-y\in\Km_v \quad (v\in F\setminus E).
\]
Then $y\in J_E$. For $v\in E\cap F$, we have
$x,y\in\Km_v$, hence $x-y\in\Km_v$.
Thus $x-y\in J_F$, proving $x\in J_E+J_F$.
The same argument, with $\Km_v$ replaced by $\CO_v$,
proves the first equality.
\end{proof}
The following result follows from Rumely's construction
\cite[Sections~2--3, pp.~203--204]{rumely1980undecidability}.
For completeness, we give a short proof using the semilocal
rings of Proposition~\ref{dittmann semilocal ring}.
\begin{lemma}\label{lem:fixed-place}
For every fixed $v\in V_K$, both $\CO_v$ and $\Km_v$ are diophantine
over $K$.
\end{lemma}
\begin{proof}
Choose distinct places $s,t\in V_K\setminus\{v\}$.
Prescribe local invariants $1/\ell$ at $v$, $-1/\ell$ at $s$,
and $0$ at every other place. Since their sum is zero,
the exactness of \eqref{eq:reciprocity} gives a class
$\beta\in\Br(K)$ with these local invariants.
By injectivity of the localization map, $\ell\beta=0$;
since the localization of $\beta$ at $v$ has order $\ell$,
the period of $\beta$ is exactly $\ell$.
By period--index equality (Lemma \ref{lem:period-index} (2)), the central division algebra $B$
representing $\beta$ has degree $\ell$.
Its ramification set is exactly $\{v,s\}$, and it splits
at every real place.
The same construction with $t$ in place of $s$ gives a
central division algebra $E$ of degree $\ell$ with
ramification set $\{v,t\}$, also split at every real place.

Proposition~\ref{dittmann semilocal ring} and
Lemma~\ref{lem:sum-intersections} give
\[
 T_B+T_E
 =(\CO_v\cap\CO_s)+(\CO_v\cap\CO_t)
 =\CO_v.
\]
Since $B$ and $E$ are fixed, $T_B$ and $T_E$ are
diophantine over $K$, and closure under sums shows that
$\CO_v$ is diophantine over $K$.

Finally, choose $\pi_v\in K$ with $v(\pi_v)=1$.
Then
\[
 \Km_v=\pi_v\CO_v,
\]
so $\Km_v$ is also diophantine over $K$.
\end{proof}\begin{lemma}\label{lem:local-splitting}
Let $w$ be a finite place of $K$ and let
$C\in\Cc_\ell^0(K)$.
Then $C\otimes_K K_w$ is split if and only if
\begin{equation}\label{eq:local-splitting}
\exists z,y\in K:\quad
z\in T_C,\quad yz=1,\quad y\in\Km_w.
\end{equation}
For fixed $w$, this condition is diophantine in the
structure constants of $C$.

If $\ell=2$, the same assertions hold for every
$C\in\Cc_2(K)$, with $S_C$ and $T_C$ defined by the
same formulas as above.
\end{lemma}

\begin{proof}
First suppose that $C$ splits at every real place.
By Proposition~3.1,
\[
T_C=\bigcap_{v\in\Delta(C)}\mathcal O_v.
\]
If $z\in T_C$ and $z^{-1}\in\Km_w$, then $w(z)<0$,
so $w\notin\Delta(C)$ and $C_w$ splits.
Conversely, if $C_w$ splits, weak approximation provides
$z\in K^\times$ with $w(z)<0$ and $v(z)\ge0$
for every $v\in\Delta(C)$.
Then $z\in T_C$ and $z^{-1}\in\Km_w$.

It remains to consider $\ell=2$ when $C$ ramifies
at a real place.
Every element of $T_C$ is still integral at every finite
ramified place, so the preceding soundness argument applies.
For the converse, suppose that $C_w$ splits.
By weak approximation, choose $z\in K^\times$ with
$w(z)<0$, integral at every finite ramified place, and
satisfying $|\sigma(z)|<1$ at every real embedding $\sigma$
where $C$ ramifies.
At each finite ramified place $v$,
\cite[Proposition~2.6(2)]{Dittmann2018} provides
$t_v\in K_v$ such that $t_v,t_v+z$ lie in the interior
of $S(C_v/K_v)$.
Such a choice is also possible at every real ramified place,
where the trace set is $[-2,2]$.
Weak approximation therefore gives $t\in K$ such that
$t,t+z$ belong to the local trace sets at all ramified places.
At split places the trace set is the whole field.
By \cite[Proposition~2.1]{Dittmann2018},
$t,t+z\in S_C$, hence $z\in T_C$.
Since $z^{-1}\in\Km_w$, this proves the converse.

Finally, membership in $T_C$ is diophantine uniformly
in the structure constants of $C$, by its definition and
Lemma~\ref{lem:reduced-functions}.
For fixed $w$, membership in $\Km_w$ is diophantine by
Lemma~\ref{lem:fixed-place}.
Thus the certificate is diophantine uniformly in the
structure constants of $C$.
\end{proof}

For a subset $I\subseteq K$, define the certificate
\begin{equation}\label{eq:certificate}
 \mathscr S(C,I):\quad
 \exists z,y\in K:\ z\in T_C,\quad yz=1,\quad y\in I.
\end{equation}
If $I\subseteq\Km_w$, this certificate implies splitting at $w$.
If $I=\Km_w$, it is equivalent to splitting at $w$. Its definability
is uniform whenever membership in $I$ is uniformly diophantine.

\section{Local witness principles}
\label{sec:witness}
\begin{proposition}[General local witness principle]\label{prop:general-witness}
Let $X\subseteq K^d$ be diophantine. Suppose a property $P(x)$ on $X$
has local conditions $P_w(x)$ indexed by a set $V$, with
\[
 P(x)\quad\Longleftrightarrow\quad
 \text{there exists }w\in V\text{ such that }P_w(x).
\]
For $1\le i\le h$, let $P_i\subseteq K^{m_i}$ and
$D_i\subseteq X\times P_i\subseteq K^{d+m_i}$ be diophantine sets.
We use $P_i$ as a set of auxiliary parameters and regard $D_i$ as
a relation between $x\in X$ and $t\in P_i$, writing $D_i(x,t)$
for $(x,t)\in D_i$.
To each $t\in P_i$, associate a nonempty subset
$E_{i,t}\subseteq V$.
Assume:
\begin{enumerate}
\item[(G1)] If $D_i(x,t)$ holds, then $P_w(x)$ holds for every
$w\in E_{i,t}$.
\item[(G2)] For every $w\in V$ and $x\in X$ satisfying $P_w(x)$,
there are an index $i$ and $t\in P_i$ with $E_{i,t}=\{w\}$
and $D_i(x,t)$.
\end{enumerate}
Then $\{x\in X:P(x)\}$ is diophantine, defined by
\[
 \bigvee_{i=1}^h \exists t\in P_i:\ D_i(x,t).
\]
\end{proposition}
\begin{proof}
Suppose $D_i(x,t)$ holds for some $i$ and $t\in P_i$.
Since $E_{i,t}$ is nonempty, (G1) gives an index $w\in V$
such that $P_w(x)$ holds, and hence $P(x)$ holds.
Conversely, if $P(x)$ holds, choose $w\in V$ such that
$P_w(x)$ holds. By (G2), there exist $i$ and $t\in P_i$
such that $D_i(x,t)$ holds.
Thus the displayed formula defines $\{x\in X:P(x)\}$.
This set is diophantine because it is a finite union of
projections of diophantine sets.
\end{proof}
\begin{remark}
The sets $E_{i,t}$ are used only to verify (G1)--(G2);
they need not be definable, and the defining formula does not
quantify over $V$.
The singleton condition $E_{i,t}=\{w\}$ in (G2) is a convenient
sufficient condition: the proof only requires that, whenever
$P(x)$ holds, there exist $i$ and $t\in P_i$ such that
$D_i(x,t)$ holds.
The substantive task in applications is to construct the
diophantine sets $P_i$ and $D_i$, together with the associated
sets $E_{i,t}$, satisfying (G1)--(G2).
\end{remark}
\begin{proposition}[Brauer detection]\label{local-diophantine-correspondence}
Let $S\subset V_K$ be finite, and let
$S_1,\ldots,S_h\subset V_K\setminus S$ be finite sets with
$\bigcap_{i=1}^h S_i=\varnothing$.
For each $1\le i\le h$, let $P_i\subseteq K^{m_i}$ be a
diophantine set of auxiliary parameters.
To each $t\in P_i$, associate an algebra
$R_{i,t}\in\Cc_\ell^0(K)$ and a subset $I_{i,t}\subseteq K$.
Assume that the sets
\[
 \{(t,\mathbf r):t\in P_i,\ 
   \mathbf r\text{ is the structure-constant table of }R_{i,t}\}
 \quad\text{and}\quad
 \{(t,z):t\in P_i,\ z\in I_{i,t}\}
\]
are diophantine over $K$.
Suppose further that:
\begin{enumerate}
\item[(W1)]
For every $1\le i\le h$ and $t\in P_i$, there exists
$w\in V_K\setminus(S\cup S_i)$ such that
\[
 w\in\Delta(R_{i,t})
 \quad\text{and}\quad
 I_{i,t}\subseteq\Km_w.
\]
\item[(W2)]
For every $1\le i\le h$ and
$w\in V_K\setminus(S\cup S_i)$, there exists $t\in P_i$
such that
\[
 w\in\Delta(R_{i,t})
 \quad\text{and}\quad
 I_{i,t}=\Km_w.
\]
\end{enumerate}
Then the set
\[
 \{Q\in\Cc_\ell^0(K):Q\text{ is not split over }K\},
\]
viewed as a set of structure-constant tables, is diophantine
over $K$.
\end{proposition}
\begin{proof}
By Lemma~\ref{lem:real-places}, the set $\Cc_\ell^0(K)$ is diophantine.
We apply Proposition~\ref{prop:general-witness} with
$X=\Cc_\ell^0(K)$, with $V$ the set of finite places of $K$,
and with $P_w(Q)$ the condition that $Q_w$ is nonsplit.

For each $v\in S$, choose a degree-$\ell$ algebra
$D_v\in\Cc_\ell^0(K)$ with $\inv_v(D_v)=1/\ell$,
as in the proof of Lemma~\ref{lem:fixed-place}.
Adjoin the pair $(D_v,\Km_v)$ as a family with a single
parameter and associated witness set $\{v\}$.
Membership in $\Km_v$ is diophantine by the same lemma.

For each original parameter $t\in P_i$, define
\[
 E_{i,t}
 =
 \{w\in V_K\setminus(S\cup S_i):
   w\in\Delta(R_{i,t}),\ I_{i,t}\subseteq\Km_w\}.
\]
This set is nonempty by (W1).
If $t$ is chosen as in (W2) for a place $w$, then
$E_{i,t}=\{w\}$.
Indeed, $w\in E_{i,t}$, and any $u\in E_{i,t}$ satisfies
$\Km_w=I_{i,t}\subseteq\Km_u$.
For distinct finite places $w$ and $u$, weak approximation
provides an element with positive $w$-valuation and negative
$u$-valuation, contradicting this inclusion.

For each family, including the adjoined singleton families,
write $(R,I)$ for the pair associated with its parameter.
Use the detector
\begin{equation}\label{eq:witness-formula}
 \bigvee_{k=1}^{\ell-1}
 \exists\, C\in\Cc_\ell^0(K):
 \qquad
 [C]=k[Q]+[R]
 \quad\text{and}\quad
 \mathscr S(C,I).
\end{equation}
Suppose this detector holds.
At every associated witness place $w$, the inclusion
$I\subseteq\Km_w$ and the certificate $\mathscr S(C,I)$
imply that $C_w$ splits.
Since $R_w$ is nonsplit, the equality
$[C_w]=k[Q_w]+[R_w]$ implies that $Q_w$ is nonsplit.
This verifies (G1).

Conversely, suppose $Q_w$ is nonsplit at a finite place $w$.
If $w\in S$, use the adjoined pair $(D_w,\Km_w)$.
Otherwise, since $\bigcap_i S_i=\varnothing$, choose $i$
with $w\notin S_i$ and apply (W2).
In either case, the resulting pair $(R,I)$ satisfies
\[
 E=\{w\},\qquad R_w\text{ is nonsplit},\qquad I=\Km_w.
\]
Both $[Q_w]$ and $[R_w]$ are nonzero elements of the cyclic
group $\Br(K_w)[\ell]$ of order $\ell$.
Hence some $1\le k<\ell$ satisfies
\[
 k[Q_w]+[R_w]=0.
\]
The global class $k[Q]+[R]$ has period dividing $\ell$.
By period--index equality over global fields, it has a
degree-$\ell$ representative $C$, using $M_\ell(K)$ when
the class is zero.
Since $Q$ and $R$ split at every real place, so does $C$;
thus $C\in\Cc_\ell^0(K)$.
Moreover, $C_w$ splits, and $I=\Km_w$, so
Lemma~\ref{lem:local-splitting} gives $\mathscr S(C,I)$.
This verifies (G2).

The hypotheses on the parameter families, the results of the
preceding section, and Lemma~7.1 show that the detectors
\eqref{eq:witness-formula} are diophantine.
Finally, for $Q\in\Cc_\ell^0(K)$, the Brauer local--global
theorem implies that $Q$ is nonsplit if and only if $Q_w$
is nonsplit at some finite place $w$.
Proposition~\ref{prop:general-witness} therefore gives the
desired conclusion.
\end{proof}

\begin{remark}\label{rem:fixed-exceptions}
The singleton families show explicitly how a fixed finite exceptional
set is handled. In particular, no definability of a local nonnorm set
is needed at places of residue characteristic $\ell$.
\end{remark}

\section{Adapted symbol presentations}
\label{sec:adapted}
\subsection{Artin--Schreier presentations}
Assume first that $\chara(K)=p>0$ and $\ell=p$.
For $a\in K$ and $b\in K^\times$, let $[a,b)_p$ denote the algebra
with generators $X,Y$ and relations
\[
 X^p-X=a,\qquad Y^p=b,\qquad YXY^{-1}=X+1.
\]
It is a central simple algebra of degree $p$, with basis
$X^iY^j$ for $0\le i,j<p$. Its structure constants are polynomial
in $a,b$, as is seen by reducing products using these relations.
\begin{definition}
A symbol presentation
\[
R=[a,b)_p,\qquad a\in K,\quad b\in K^\times,
\]
is \emph{adapted} if
\[
a\in T_R.
\]
Equivalently, by Proposition~\ref{dittmann semilocal ring},
$a$ is integral at every place in $\Delta(R)$.
\end{definition}

The set of adapted pairs
\[
\{(a,b)\in K\times K^\times:a\in T_{[a,b)_p}\}
\]
is diophantine. Indeed, the condition $b\neq0$ is diophantine,
the structure constants of $[a,b)_p$ are polynomial in $a$ and $b$,
and membership in $T_R$ is uniformly diophantine by
Proposition~\ref{dittmann semilocal ring}.
\begin{lemma}\label{lem:adapted-existence}
Every degree-$p$ central simple algebra over $K$ has an adapted symbol
presentation.
\end{lemma}
\begin{proof}
For the split algebra use $[0,1)_p$. Otherwise $\Delta(A)$ is
nonempty. At each $v\in\Delta(A)$ choose a residue class $\bar a_v$
outside the image of $u\mapsto u^p-u$ on the finite residue field.
This is possible since the map has kernel $\F_p$ and hence is not
surjective. By weak approximation choose $a\in K$ integral at these
places with the prescribed residues. An Artin--Schreier polynomial of
prime degree over a finite field is either split or irreducible.
Thus $Z^p-Z-a$ defines the unramified degree-$p$ extension at each
$v\in\Delta(A)$, and is irreducible over $K$ as well.

Put $L=K[Z]/(Z^p-Z-a)$. It is cyclic of degree $p$, with generator
$\sigma(\bar Z)=\bar Z+1$. At every ramified place its local degree
is $p$, so restriction multiplies the local invariant of $A$ by $p$
and kills it. At all other places $A$ already splits. Brauer
reciprocity over $L$ therefore implies that $A\otimes_K L$ is split.
Since $[L:K]=\deg(A)=p$, the maximal-subfield embedding criterion
and the cyclic algebra construction give
$A\simeq(L/K,\sigma,b)=[a,b)_p$ for some $b\ne0$.
Finally $a$ is integral at every place of $\Delta(A)$, so $a\in T_A$.
\end{proof}

\begin{lemma}\label{adapted symbol local behaviour}
If $R=[a,b)_p$ is adapted and $v\in\Delta(R)$, then
$Z^p-Z-a$ defines the unramified degree-$p$ extension of $K_v$, and
$v(b)\not\equiv0\pmod p$.
\end{lemma}
\begin{proof}
We have $a\in\CO_v$. If the residue polynomial had a root, Hensel's
lemma would give a root in $K_v$, and the symbol would split there.
Its residue polynomial is therefore irreducible, giving the stated
unramified extension. The norm group of an unramified extension of
local fields of degree $p$ consists of the elements of valuation
divisible by $p$; see \cite{SerreLocalFields}. Since a cyclic algebra splits precisely when its
parameter is a norm, $v(b)$ is not divisible by $p$.
\end{proof}

\subsection{Kummer presentations}
Assume $\ell\ne\chara(K)$ and $\mu_\ell\subset K$, and fix a primitive
$\ell$th root of unity $\zeta$. Let $S$ be the set of finite places of
residue characteristic $\ell$; it is empty for function fields. Put
\[
 (a,b)_\zeta=K\langle X,Y\mid X^\ell=a,
 Y^\ell=b,\ YX=\zeta XY\rangle,\qquad a,b\in K^\times.
\]
This is a central simple algebra of degree $\ell$, even when $a$ is
an $\ell$th power. Its table in the monomial basis is polynomial in
$a,b$. Call this presentation \emph{adapted} if $R=(a,b)_\zeta$ splits
at all real places and at every place of $S$, and
\[
 a,a^{-1}\in T_R.
\]
These conditions are diophantine by Lemmas~\ref{lem:real-places}
and~\ref{lem:local-splitting}.

\begin{lemma}\label{lem:kummer-adapted}
Every degree-$\ell$ algebra split at all real places and at $S$ has
an adapted Kummer presentation. For every such presentation and
$v\in\Delta(R)$, the extension $K_v(\sqrt[\ell]{a})/K_v$ is
unramified of degree $\ell$, and $\ell\nmid v(b)$.
\end{lemma}
\begin{proof}
For a split algebra use $(1,1)_\zeta$. Otherwise, for every
$v\in\Delta(A)$, the residue characteristic is different from $\ell$
and the residue field contains $\mu_\ell$. Its multiplicative group
therefore has non-$\ell$th powers. By weak approximation choose $a$
which is a unit at these places with non-$\ell$th-power residue.
The polynomial $Z^\ell-a$ defines the unramified degree-$\ell$
extension at every such place; in particular it is irreducible over
$K$. Here irreducibility follows directly from local Kummer theory;
we do not use a diophantine irreducibility criterion.

The cyclic extension $L=K(\sqrt[\ell]{a})$ splits $A$: its local
degree is $\ell$ wherever $A$ ramifies, and $A$ splits at all other
places. The embedding criterion for maximal subfields and the cyclic
algebra construction give $A\simeq(a,b)_\zeta$ for some $b\ne0$,
using the generator sending $\sqrt[\ell]{a}$ to
$\zeta\sqrt[\ell]{a}$. The chosen $a$ is a unit at every ramified
place, so the presentation is adapted.

Conversely, in an adapted presentation $a$ is a unit at every
$v\in\Delta(R)$. If its residue were an $\ell$th power, Hensel's
lemma would give $a\in K_v^{\times\ell}$ and split $R_v$.
Thus the local Kummer extension is unramified of degree $\ell$.
Its norm group consists of elements of $\ell$-divisible valuation,
so nonsplitting forces $\ell\nmid v(b)$; see
\cite{SerreLocalFields,GilleSzamuely2017}.
\end{proof}

\subsection{The common input}
In the Artin--Schreier case take $S=\varnothing$; in the Kummer case
take $S$ as above. In both cases the allowed symbol presentations
form a diophantine parameter set, give tables in $\Cc_\ell^0(K)$
split at $S$, and satisfy
\begin{equation}\label{eq:common-valuation}
 \ell\nmid v(b)\qquad(v\in\Delta(R)).
\end{equation}
Every degree-$\ell$ algebra split at $S$ and all real places has such
a presentation. These are the only properties of the two symbol
presentations needed below.

\section{A common radical formula and ramification matching}
\label{sec:synchronization}

\subsection{A uniform radical formula}
For either adapted presentation from the preceding section, put
\[
 J_R=\bigcap_{v\in\Delta(R)}\Km_v,\qquad
 T_R^\times=\{u\in K^\times:u,u^{-1}\in T_R\},
\]
and, for $1\le k<\ell$, define
\begin{equation}\label{eq:Ck}
 C_k=\{0\}\cup
 \bigl(b^kK^{\times \ell}T_R^\times
             \cap(1-K^{\times \ell}T_R^\times)\bigr).
\end{equation}
These sets are uniformly diophantine in $(a,b)$. Their valuation
properties use only \eqref{eq:common-valuation}.

\begin{lemma}\label{valuation description Ck}
Suppose $\Delta(R)\ne\varnothing$. For $x\in K^\times$, membership
in $C_k$ is equivalent to
\begin{equation}\label{eq:Ck-valuations}
 v(x)>0,\qquad v(x)\equiv k v(b)\pmod \ell
 \quad\text{for every }v\in\Delta(R).
\end{equation}
\end{lemma}\begin{proof}
Suppose that $x\in C_k$ is nonzero. By \eqref{eq:Ck}, we may write
\[
x=b^k z^\ell u=1-(z')^\ell u'
\]
with $z,z'\in K^\times$ and $u,u'\in T_R^\times$.
Fix $v\in\Delta(R)$. Since $v(u)=v(u')=0$, the first expression gives
\[
v(x)\equiv k v(b)\not\equiv 0\pmod{\ell}.
\]
If $v(x)<0$, then $v(1-x)=v(x)$, while the second expression gives
$v(1-x)=\ell v(z')$, a contradiction. Thus $v(x)>0$.

Conversely, suppose that $v(x)>0$ and
$v(x)\equiv k v(b)\pmod{\ell}$ for every $v\in\Delta(R)$.
By weak approximation, choose $z\in K^\times$ such that
\[
v(z)=\frac{v(x)-k v(b)}{\ell}
\qquad\text{for every }v\in\Delta(R).
\]
Then $u=x/(b^k z^\ell)$ has valuation zero at every such place,
so $u\in T_R^\times$. Moreover, $v(1-x)=0$ at every such place
because $v(x)>0$. Since $\Delta(R)$ is nonempty, this also gives
$1-x\ne0$, hence $1-x\in T_R^\times$. Taking $z'=1$ and $u'=1-x$
gives the second expression in \eqref{eq:Ck}. Therefore $x\in C_k$.
\end{proof}

\begin{lemma}\label{odd radical lemma}\label{lem:radical}
For every such adapted presentation, including $\ell=2$, one has
\begin{equation}\label{eq:radical}
 J_R=\sum_{k=1}^{\ell-1}(C_k+C_k).
\end{equation}
Consequently, membership in $J_R$ is uniformly diophantine in $(a,b)$.
\end{lemma}
\begin{proof}
First assume $\Delta(R)\ne\varnothing$. The right-hand side is
contained in $J_R$ by Lemma~\ref{valuation description Ck}.
Let $x\in J_R$. We construct, in each completion at a place
$v\in\Delta(R)$, nonzero elements $c_{k,1,v},c_{k,2,v}$ whose sum
over $k$ is $x$ and whose positive valuations are congruent to
$kv(b)$ modulo $\ell$.

Let $r_{k,v}\in\{1,\ldots,\ell-1\}$ be the least positive residue of
$kv(b)$ modulo $\ell$. There is a unique $k_v$ with $r_{k_v,v}=1$.
For every $k\ne k_v$, choose a pair $c,-c$ with $v(c)=r_{k,v}$.
For the remaining pair, if $x=0$ choose $c,-c$ with $v(c)=1$.
If $x\ne0$ and $v(x)>1$, choose $c$ of valuation $1$ and use
$(c,x-c)$; both entries have valuation $1$. If $v(x)=1$, choose
$c$ of valuation $\ell+1$ and use $(c,x-c)$; the valuations are
$\ell+1$ and $1$. These constructions work also in residue
characteristic two.

Choose all but one of the $2(\ell-1)$ global entries in $K$ sufficiently
close to their local choices at every $v\in\Delta(R)$, using weak
approximation. Define the final entry to be $x$ minus the sum of the
others. It too is sufficiently close to its prescribed nonzero local
value, since the local entries sum to $x$. Thus all global entries
have the required positive valuations and congruences. By
Lemma~\ref{valuation description Ck}, each belongs to its respective
$C_k$, proving \eqref{eq:radical}.

If $\Delta(R)=\varnothing$, then $T_R=J_R=K$ and
$C_k=K\setminus\{1\}$ for every $k$. Since $K$ is infinite,
$(K\setminus\{1\})+(K\setminus\{1\})=K$, so the formula remains
valid. The definability assertion follows from \eqref{eq:Ck}.
\end{proof}

\subsection{Ramification matching via a fixed Brauer class}
\label{subsec:ramification-matching}
Fix distinct places $s,t\in V_K\setminus S$. By \eqref{eq:reciprocity}, choose a
fixed degree-$\ell$ algebra $D=D_{s,t}$ with
\[
 \inv_s(D)=1/\ell,\qquad \inv_t(D)=-1/\ell,
 \qquad \inv_v(D)=0\quad(v\ne s,t).
\]
Write $R(a,b)$ for the appropriate symbol algebra.
Let $P_{s,t}\subseteq K^4$ consist of quadruples
$\mathbf t=(a_1,b_1,a_2,b_2)$ such that
\begin{equation}\label{eq:synchronization-parameters}
 \begin{gathered}
 b_1b_2\ne0,\qquad R_1=R(a_1,b_1),\quad R_2=R(a_2,b_2),\\
 (a_1,b_1),(a_2,b_2)\text{ are adapted},\qquad
 [R_1]-[R_2]=[D],\\
 R_1\otimes_K K_t\text{ is split},\qquad
 R_2\otimes_K K_s\text{ is split}.
 \end{gathered}
\end{equation}
This set is diophantine: the class equality can be expressed as
$[R_1]=[R_2\otimes_K D]$, and the two fixed-place conditions are
covered by Lemma~\ref{lem:local-splitting}. For $\mathbf t\in P_{s,t}$,
put
\[
 R_{\mathbf t}=R_1,\qquad I_{\mathbf t}=J_{R_1}+J_{R_2}.
\]
Membership in $I_{\mathbf t}$ is uniformly diophantine by
Lemma~\ref{lem:radical}.

\begin{lemma}[Ramification matching]\label{odd synchronization}
The family just defined satisfies (W1) and (W2) of
Proposition~\ref{local-diophantine-correspondence} outside $S$, with seed set $\{s,t\}$.
\end{lemma}
\begin{proof}
For any parameter in $P_{s,t}$, the class equality and the two
splitting conditions give
\[
 \inv_s(R_1)=1/\ell,\quad \inv_s(R_2)=0,\qquad
 \inv_t(R_1)=0,\quad \inv_t(R_2)=1/\ell.
\]
At every $v\ne s,t$, the invariants of $R_1$ and $R_2$ are equal.
Consequently, there is a finite set $E\subseteq V_K\setminus(S\cup\{s,t\})$
such that
\[
 \Delta(R_1)=\{s\}\cup E,\qquad
 \Delta(R_2)=\{t\}\cup E.
\]
This set $E$ is nonempty: otherwise $R_1$ would have total invariant
$1/\ell$, contradicting reciprocity. By
Lemma~\ref{lem:sum-intersections},
\begin{equation}\label{eq:synchronized-ideal}
 I_{\mathbf t}=\bigcap_{v\in E}\Km_v.
\end{equation}
Any $w\in E$ proves (W1).

For (W2), let $w\notin S\cup\{s,t\}$. Prescribe a degree-$\ell$ algebra $B_1$ with
invariants $1/\ell$ at $s$, $-1/\ell$ at $w$, and zero elsewhere; similarly
prescribe $B_2$ with invariants $1/\ell$ at $t$, $-1/\ell$ at $w$, and zero
elsewhere. Reciprocity and period--index equality give these algebras,
and injectivity in \eqref{eq:reciprocity} gives
$[B_1]-[B_2]=[D]$. These algebras split at $S$ and all real places. By
Lemmas~\ref{lem:adapted-existence} and~\ref{lem:kummer-adapted}, as
appropriate, choose
adapted presentations of both algebras. They define a parameter in
$P_{s,t}$ for which $E=\{w\}$, so
$I_{\mathbf t}=\Km_w$ and $w\in\Delta(R_{\mathbf t})$.
\end{proof}

\begin{remark}
The fixed Brauer-class difference in
\eqref{eq:synchronization-parameters} forces the ramification sets
of the two auxiliary algebras to agree outside the prescribed
seed places. We refer to this construction as ramification matching.
The nonemptiness of the common ramification set follows from Brauer
reciprocity, rather than from an additional nonvanishing condition
on the variable algebras.
\end{remark}

\section{Completion of the main theorem}
\label{sec:completion}

\subsection{Real places}
\begin{lemma}\label{lem:real-places}
For a global field $K$ and prime $\ell$, the set $\Cc_\ell^0(K)$ is
diophantine. For each fixed real place $v$, the set of degree-two
tables nonsplit at $v$ is also diophantine.
\end{lemma}
\begin{proof}
Only degree two over a number field requires proof. For each real
embedding $v$, the set $P_v=\{x\in K^\times:v(x)>0\}$ is
diophantine by \cite[Lemma~4.1]{KTU}. This lemma uses reduced norm
values of a fixed quaternion algebra, not nonnorm definability.
Every degree-two central simple algebra has a quaternion presentation
$A\simeq(a,b)_{-1}$, $a,b\ne0$. At $v$ it splits if and only if
$a\in P_v$ or $b\in P_v$, and is nonsplit if and only if
$-a,-b\in P_v$. Quantify a presentation and an isomorphism as in
Lemma~\ref{lem:algebra-isomorphism}. Taking the conjunction of the
splitting conditions over the finitely many real places proves the
first assertion; the second follows from the nonsplitting condition.
\end{proof}

\subsection{Prime degrees with symbol presentations}
\begin{proposition}\label{prop:symbol-prime}
Suppose either $\ell=\chara(K)$, or $\ell\ne\chara(K)$ and
$\mu_\ell\subset K$. Then $\Nc_\ell(K)$ is diophantine.
\end{proposition}
\begin{proof}
Choose four distinct finite places $s_1,t_1,s_2,t_2$ outside $S$.
For $i=1,2$, take the matched ramification family
$P_{s_i,t_i}$ defined by \eqref{eq:synchronization-parameters},
with seed set $\{s_i,t_i\}$. The parameter conditions are
diophantine, and Lemma~\ref{lem:radical} gives a uniform
diophantine definition of the associated ideals.

By Lemma~\ref{odd synchronization}, each family satisfies
\textup{(W1)} and \textup{(W2)} outside $S$ with its respective
seed set. Since the seed sets are disjoint, every finite place
$w\notin S$ lies outside at least one of them and is therefore
covered by \textup{(W2)} for the corresponding family.
Proposition~\ref{local-diophantine-correspondence} now gives a
diophantine definition of the nonsplit locus in $\Cc_\ell^0(K)$,
including algebras that are nonsplit at a place in $S$.

If $\ell=2$ and $K$ has real places, taking the union with the
finitely many real nonsplitting loci from Lemma~\ref{lem:real-places}
gives $\Nc_\ell(K)$. In all other cases
$\Cc_\ell^0(K)=\Cc_\ell(K)$, so the result follows.
\end{proof}

\begin{remark}
For $\ell=2$ the common formulas are $T_A=S_A-S_A$ and
$J_R=C_1+C_1$. The radical argument allows valuations $1$ and $3$
when the target has valuation $1$, and therefore also works over
residue field $\F_2$.
\end{remark}

\subsection{Removing roots of unity}
\begin{lemma}\label{lem:descent}
Let $L/K$ be a fixed finite extension of degree prime to $\ell$.
If $\Nc_\ell(L)$ is diophantine over $L$, then $\Nc_\ell(K)$ is
diophantine over $K$.
\end{lemma}
\begin{proof}
Restriction on $\ell$-torsion Brauer classes is injective: the
composition of restriction and corestriction is multiplication by
$[L:K]$, which is invertible modulo $\ell$. Hence
\[
 A\text{ is nonsplit over }K
 \quad\Longleftrightarrow\quad A\otimes_K L\text{ is nonsplit over }L.
\]
Choose a fixed $K$-basis of $L$. An existential polynomial formula
over $L$, with parameters in $L$, becomes a finite system of
polynomial equations over $K$ by expanding every quantified variable
in this basis and equating coordinates. Apply this to the formula for
$\Nc_\ell(L)$ and the scalar-extended structure constants of $A$.
\end{proof}

\begin{proof}[Proof of Theorem~\ref{prime version}]
The characteristic-degree case follows from
Proposition~\ref{prop:symbol-prime}. Otherwise set $L=K(\mu_\ell)$.
The degree $[L:K]$ divides $\ell-1$, so is prime to $\ell$.
Apply Proposition~\ref{prop:symbol-prime} over $L$ and then
Lemma~\ref{lem:descent}.
\end{proof}

\subsection{Reduction to prime degree}
\begin{proposition}\label{prop:prime-reduction}
Over a global field $K$, Theorem~\ref{prime version} implies
Theorem~\ref{main theorem}.
\end{proposition}
\begin{proof}
For $A$ of degree $n$, we claim that $A$ is nonsplit if and only if
there are a prime $\ell\mid n$, an integer $1\le m\le n$, and a
nonsplit algebra $B$ of degree $\ell$ such that $[B]=m[A]$.
If $d=\per(A)>1$, choose $\ell\mid d$ and put $m=d/\ell$.
The class $m[A]$ has period and index $\ell$, so its division algebra
is such a $B$. Conversely, $[B]=m[A]\ne0$ implies $[A]\ne0$.
For fixed $\ell,m$, the condition is existential: quantify the table
of $B$, impose membership in $\Nc_\ell(K)$, and use
Lemma~\ref{lem:algebra-isomorphism} to express $[B]=[A^{\otimes m}]$.
There are only finitely many pairs $(\ell,m)$ to consider.
\end{proof}

\begin{proof}[Proof of Theorem~\ref{main theorem}]
Combine Theorem~\ref{prime version} with
Proposition~\ref{prop:prime-reduction}.
\end{proof}

\section{Applications and alternative proofs}
\label{sec:applications}
The following applications illustrate two uses of the construction.
Nonintegrality uses the localization families directly, with a
valuation detector in place of Brauer detection. Nonnorms, quadratic forms, and the extension certificate for rootless
polynomials use nonsplitting certificates. General four-dimensional
forms also require a local quadratic splitting condition. The earlier results are attributed below;
the claims of simplification concern the indicated proof steps, not
quantifier counts or the removal of all classical arithmetic inputs.

\subsection{The complement of the integers}
\begin{corollary}[Koenigsmann]\label{cor:nonintegers}
The set $\Q\setminus\Z$ is diophantine over $\Q$.
\end{corollary}
\begin{proof}
Apply the Kummer ramification-matching construction over $\Q$ with
$\ell=2$ and $S=\{2\}$, as in Lemma~\ref{odd synchronization}.
Choose two pairs of odd primes with disjoint supports. Their matched ramification
families have uniformly diophantine ideals
\[
 I_{i,t}=\bigcap_{v\in E_{i,t}}\Km_v,
 \qquad \varnothing\ne E_{i,t}\subset V_\Q\setminus\{2\},
\]
and can isolate every odd prime by a suitable choice of $i,t$.
For these families use the detector
\[
 D_i(x,t):\quad \exists y\in\Q:\ xy=1,\quad y\in I_{i,t}.
\]
It forces $v(x)<0$ at every $v\in E_{i,t}$, and conversely holds
whenever $v(x)<0$ and $I_{i,t}=\Km_v$. Add the fixed ideal
$\Km_2$, which is diophantine by Lemma~\ref{lem:fixed-place},
to detect a pole at $2$. Thus
\[
 x\in\Q\setminus\Z\quad\Longleftrightarrow\quad
 \bigl(\exists y:\ xy=1,\ y\in\Km_2\bigr)
 \ \lor\ \bigvee_i\exists t\in P_i:\ D_i(x,t).
\]
This is an instance of Proposition~\ref{prop:general-witness} with
local condition $v(x)<0$. In particular the certificate excludes
$x=0$, as required.
\end{proof}

\begin{remark}
The diophantineness of the complement of the ring of
integers over number fields was established by
Koenigsmann for $\Q$ \cite{Koenigsmann2016}
and extended to arbitrary number fields by
Park \cite{Park2013}. Related results for
global function fields and rings of $S$-integers
were obtained in \cite{KTU,Daans2021}.

The argument above provides a unified alternative
proof over arbitrary global fields. Its main
advantage is that the ramification-matching
construction treats all eligible finite places
uniformly, without requiring explicit congruence
conditions or separate calculations of local
Hilbert symbols. The underlying arithmetic input
is supplied instead by Brauer--Hasse--Noether
reciprocity and period--index equality.
\end{remark}

\begin{corollary}[Nonintegrality over global fields]
\label{cor:global-nonintegrality}
Let $K$ be a global field, and let
$S\subseteq V_K$ be finite. Then
\[
K\setminus\mathcal O_{K,S}
\]
is diophantine over $K$, where
\[
\mathcal O_{K,S}
=\bigcap_{v\in V_K\setminus S}\mathcal O_v.
\]
In particular, for every number field $K$,
the complement $K\setminus\mathcal O_K$
is diophantine over $K$.
\end{corollary}

\begin{proof}
We apply the ramification-matching construction
of Subsection~\ref{subsec:ramification-matching}
with $\ell=2$, using the adapted Kummer
presentations of Lemma~\ref{lem:kummer-adapted}
when $\chara(K)\ne2$, and the adapted
Artin--Schreier presentations of
Lemma~\ref{lem:adapted-existence} otherwise.

Let $S_0$ be the exceptional set of finite
places in this construction. Thus $S_0$
consists of the places of residue
characteristic $2$ when $\chara(K)\ne2$,
and is empty when $\chara(K)=2$.
Put $F=S\cup S_0$.

Choose four distinct finite places
$s_1,t_1,s_2,t_2$ outside $F$.
For each pair $(s_i,t_i)$, consider the
ramification-matching family constructed
in Subsection~\ref{subsec:ramification-matching},
with the additional requirement that both
auxiliary algebras split at every place
of $S\setminus S_0$.
Denote its parameter set by $P_i$.

These additional splitting conditions are
diophantine by Lemma~\ref{lem:local-splitting}.
They do not affect the existence of
parameters isolating a prescribed place
$w\notin F$: the auxiliary algebras
constructed in the proof of
Lemma~\ref{odd synchronization} are
ramified only at $\{s_i,w\}$ and
$\{t_i,w\}$, respectively, and therefore
split at every place of $S$.

Consequently, the proof of
Lemma~\ref{odd synchronization} applies
to these restricted families.
For every $\mathbf t\in P_i$, the associated
ideal is uniformly diophantine and has
the form
\[
I_{i,\mathbf t}
=\bigcap_{v\in E_{i,\mathbf t}}\Km_v,
\qquad
\varnothing\ne E_{i,\mathbf t}
\subseteq V_K\setminus F.
\]
Moreover, since the two seed pairs are
disjoint, for every $w\notin F$ there
exist $i\in\{1,2\}$ and
$\mathbf t\in P_i$ such that
\[
E_{i,\mathbf t}=\{w\},
\qquad
I_{i,\mathbf t}=\Km_w.
\]

For each family, consider the condition
\[
D_i(x,\mathbf t):
\quad
\exists y\in K:
\quad xy=1,\qquad y\in I_{i,\mathbf t}.
\]
This is diophantine uniformly in
$(x,\mathbf t)$.
If $D_i(x,\mathbf t)$ holds, then
\[
v(x)<0
\qquad(v\in E_{i,\mathbf t}).
\]
Since $E_{i,\mathbf t}$ is nonempty and
disjoint from $S$, this implies
$x\notin\mathcal O_{K,S}$.

Conversely, suppose $v_w(x)<0$ for
some $w\notin F$. Choose $i$ and
$\mathbf t\in P_i$ such that
$I_{i,\mathbf t}=\Km_w$.
Then $y=x^{-1}\in\Km_w$, so
$D_i(x,\mathbf t)$ holds.

It remains to handle the finitely many
places in $S_0\setminus S$.
For each such place $v$, the ideal
$\Km_v$ is diophantine by
Lemma~\ref{lem:fixed-place}, and
\[
v(x)<0
\quad\Longleftrightarrow\quad
\exists y\in K:
\quad xy=1,\qquad y\in\Km_v.
\]

We therefore obtain
\[
\begin{aligned}
x\in K\setminus\mathcal O_{K,S}
\quad\Longleftrightarrow\quad&
\bigvee_{v\in S_0\setminus S}
\left(
\exists y\in K:
xy=1,\ y\in\Km_v
\right)
\\
&\lor
\bigvee_{i=1}^{2}
\exists\mathbf t\in P_i:
D_i(x,\mathbf t).
\end{aligned}
\]
All conditions on the right are
diophantine, and the disjunctions
are finite. Hence
$K\setminus\mathcal O_{K,S}$
is diophantine over $K$.
\end{proof}

\subsection{Prime-degree nonnorms}
\begin{corollary}[Prime-degree nonnorms]
\label{cor:nonnorms}
Let $K$ be a global field.

If $\ell\ne\chara(K)$ is prime and
$\mu_\ell\subset K$, then the set
\[
\left\{(a,b)\in K^\times\times K^\times:
b\notin
\Norm_{K(\sqrt[\ell]{a})/K}
\bigl(K(\sqrt[\ell]{a})^\times\bigr)
\right\}
\]
is diophantine over $K$.

If $\chara(K)=p>0$, then the set
\[
\left\{(a,b)\in K\times K^\times:
\begin{array}{l}
Z^p-Z-a\text{ is irreducible over }K,\\
b\notin\Norm_{L_a/K}(L_a^\times)
\end{array}
\right\},
\qquad
L_a=K[Z]/(Z^p-Z-a),
\]
is diophantine over $K$ as well.
\end{corollary}

\begin{proof}
Fix a primitive $\ell$th root of unity
$\zeta\in K$ in the first case.
By the cyclic algebra splitting criterion
\cite[Section~4.7]{GilleSzamuely2017},
\[
(a,b)_\zeta\text{ is nonsplit}
\quad\Longleftrightarrow\quad
b\notin
\Norm_{K(\sqrt[\ell]{a})/K}
\bigl(K(\sqrt[\ell]{a})^\times\bigr).
\]
If $a$ is an $\ell$th power, both conditions
are false, so no additional restriction
on $a$ is needed.

In characteristic $p$, the Artin--Schreier
polynomial $Z^p-Z-a$ is either irreducible
or splits completely over $K$.
In the latter case $[a,b)_p$ is split.
Otherwise, $L_a/K$ is cyclic of degree $p$,
and the same splitting criterion gives
\[
[a,b)_p\text{ is nonsplit}
\quad\Longleftrightarrow\quad
b\notin\Norm_{L_a/K}(L_a^\times).
\]

In both cases, the structure constants
of the corresponding symbol algebra
are polynomial in $a$ and $b$.
The conclusion therefore follows from
Theorem~\ref{prime version}.
\end{proof}

\begin{remark}
For odd $\ell$, the first assertion is Morrison's theorem
\cite{Morrison2019}; the quadratic case in characteristic
different from two follows from earlier work, including
\cite{KTU}. The arithmetic ingredients of our construction
are closely related to those used in these earlier proofs,
particularly the use of local obstructions and diophantine
descriptions of semilocal rings and their Jacobson radicals.

The main difference lies in the organization of the argument.
Our ramification-matching construction separates the
localization mechanism from the detection of nonsplitting,
allowing the same method to treat arbitrary prime-degree
central simple algebras, including the characteristic-degree
case. Thus the corollary recovers the known nonnorm results
as consequences of a more general diophantine definability
theorem, rather than through a separate analysis of norm
equations.
\end{remark}
\subsection{Univariate polynomials without roots}
We next give an alternative derivation of Dittmann's rootlessness
theorem \cite[Theorem~5.2]{Dittmann2018}. We retain the classical
existence input used in his proof, but replace the auxiliary ideal
and ray-class construction of his Section~3 and the ensuing detector
by a relative Brauer-class certificate.

\begin{lemma}[Relative Brauer witnesses]\label{lem:relative-brauer}
Let $L/K$ be a nontrivial finite extension of global fields, of degree
$n$. There are a prime $\ell\mid n$ and a nonsplit degree-$\ell$
central simple $K$-algebra $A$ such that $A\otimes_K L$ splits.
\end{lemma}
\begin{proof}
By \cite[Corollary~4]{FeinKantorSchacher1981}, the relative Brauer
group $\Br(L/K)=\ker(\Br(K)\to\Br(L))$ is infinite. Since it is
torsion, it contains a class $\beta$ of prime order $\ell$.
Let $A$ be the central division $K$-algebra representing $\beta$.
By period--index equality over $K$, we have $\deg A=\ell$.
As $\beta\in\Br(L/K)$, the extension $L$ splits $A$, so
$\ell=\ind(A)$ divides $[L:K]=n$.
\end{proof}
\begin{corollary}[Dittmann]\label{cor:rootlessness}
For every global field $K$ and fixed $n\ge1$, the coefficient tuples
of monic degree-$n$ polynomials without a root in $K$ form a
diophantine set.
\end{corollary}
\begin{proof}
For a finite-dimensional commutative unital $K$-algebra $E$ and
$n>1$, consider
\begin{equation}\label{eq:extension-detector}
 \Psi_n(E):\quad
 \bigvee_{\substack{\ell\mid n\\\ell\text{ prime}}}
 \exists A\in\Nc_\ell(K):\quad
 A\otimes_K E\simeq M_\ell(E)\quad\text{as }E\text{-algebras}.
\end{equation}
For fixed $\dim_K E$, this is a diophantine condition in its table.
Indeed, membership in $\Nc_\ell(K)$ is given by
Theorem~\ref{prime version}. An $E$-linear isomorphism is expressed
by its matrix and inverse, together with the multiplication and
identity equations. Expressing its entries in a $K$-basis of $E$
turns these into finitely many polynomial equations over $K$.
This encoding does not require $E$ to be a field or reduced.

By Lemma~\ref{lem:relative-brauer}, $\Psi_n(L)$ holds for every
degree-$n$ field extension $L/K$. On the other hand, it cannot hold
for any $E$ admitting a $K$-algebra map $E\to K$: base change along
that map would split the very algebra required to lie in
$\Nc_\ell(K)$.

Now set $E_f=K[T]/(f)$ for a monic degree-$n$ polynomial. Its
multiplication table is polynomial in the coefficients of $f$.
A $K$-root of $f$ gives $E_f\to K$, so $\Psi_n(E_f)$ implies that
$f$ has no root. Conversely, if $f$ is irreducible and $n>1$,
then $\Psi_n(E_f)$ holds by the lemma.

Define the rootlessness conditions inductively. For $n=1$ use the
empty set. For $n>1$, use the disjunction of $\Psi_n(E_f)$ with
\[
 \bigvee_{r=2}^{n-2}\exists\text{ monic }g,h:\quad
 f=gh,\quad \deg g=r,\quad\deg h=n-r,\quad
 g,h\text{ have no root in }K.
\]
The union is empty when $n<4$. Every reducible rootless monic
polynomial factors into two rootless monic polynomials of degrees
at least two, whereas a polynomial with a root cannot so factor.
By induction all these are diophantine conditions. This is the
factorization induction of \cite[Lemma~5.1]{Dittmann2018}, with
\eqref{eq:extension-detector} replacing his extension detector.
\end{proof}

\begin{remark}
This proof illustrates how nonsplitting definability yields a
direct diophantine certificate for proper finite extensions,
and hence an alternative derivation of Dittmann's rootlessness
theorem. The certificate replaces his auxiliary ideal and
ray-class construction, while retaining the
Fein--Kantor--Schacher existence input and the factorization
induction. The simplification is therefore relative to the
availability of our nonsplitting theorem.
\end{remark}
\subsection{General quadratic forms}
\label{subsec:general-quadratics}

We use $\perp$ to denote orthogonal sum.
In characteristic different from two, we write
\[
\langle u_1,\ldots,u_m\rangle
=\sum_{i=1}^m u_iX_i^2,
\]
whereas in characteristic two, we put
\[
[1,r]=X^2+XY+rY^2.
\]
Over a field $F$ of the respective characteristic, the binary
forms $\langle1,-r\rangle$ and $[1,r]$ are the norm forms
of the quadratic etale algebras $F[T]/(T^2-r)$ and
$F[T]/(T^2+T+r)$, respectively.

We also recall the quaternion-algebra notation.
If $\chara(F)\ne2$, then for $u,v\in F^\times$,
the symbol $(u,v)_{-1}$ denotes the $F$-algebra generated
by $i,j$ subject to
\[
i^2=u,\qquad j^2=v,\qquad ji=-ij.
\]
If $\chara(F)=2$, then for $u\in F$ and $v\in F^\times$,
the symbol $[u,v)_2$ denotes the $F$-algebra generated
by $i,j$ subject to
\[
i^2+i=u,\qquad j^2=v,\qquad ji=(i+1)j.
\]
Both are central simple $F$-algebras of degree two.
The cyclic algebra splitting criterion gives
\[
\begin{aligned}
(u,v)_{-1}\text{ splits}
&\iff v\in N_{L/F}(L^\times),
&&L=F[T]/(T^2-u),\\
[u,v)_2\text{ splits}
&\iff v\in N_{L/F}(L^\times),
&&L=F[T]/(T^2+T+u),
\end{aligned}
\]
in the respective characteristics.

\begin{lemma}[Local four-dimensional criterion]
\label{lem:local-four}
Let $F$ be a nonarchimedean local field.
In characteristic different from two, put
\[
q=a\langle1,-r\rangle\perp c\langle1,-s\rangle,
\qquad
\delta=rs,\qquad B=(r,-c/a)_{-1},
\]
where $a,c,r,s\in F^\times$.
In characteristic two, put
\[
q=a[1,r]\perp c[1,s],
\qquad
\delta=r+s,\qquad B=[r,c/a)_2,
\]
where $a,c\in F^\times$ and $r,s\in F$.
Then $q$ is anisotropic if and only if $B$ is nonsplit and
$\delta\in F^{\times2}$ in the first case, or
$\delta\in\wp(F)$ in the second, where
$\wp(F)=\{z^2+z:z\in F\}$.
\end{lemma}

\begin{proof}
Suppose first that $\chara(F)\ne2$.
We write $(u,v)_F\in\{\pm1\}$ for the Hilbert symbol,
to distinguish it from the quaternion algebra $(u,v)_{-1}$.
By \cite[Theorem~2.3(3)]{hu2026}, applied to
\[
(a_1,a_2,a_3,a_4)=(a,-ar,c,-cs),
\]
the form $q$ is anisotropic if and only if
\[
rs\in F^{\times2}
\qquad\text{and}\qquad
(a,-ar)_F=-(-c,cs)_F.
\]
Under the square condition, $r$ and $s$ have the same
square class, and
\[
\begin{aligned}
(a,-ar)_F\,(-c,cs)_F
  &=(a,r)_F\,(-c,s)_F\\
  &=(a,r)_F\,(-c,r)_F\\
  &=(r,-ac)_F\\
  &=(r,-c/a)_F.
\end{aligned}
\]
Thus the second condition is precisely nonsplitting of $B$.

Suppose now that $\chara(F)=2$.
For $t\in F$, put
\[
L_t=F[T]/(T^2+T+t),
\qquad
H_t=N_{L_t/F}(L_t^\times).
\]
Since $[1,t]$ is the norm form of $L_t$, we have
\[
q\text{ is isotropic}\iff c/a\in H_rH_s.
\]
Indeed, if both algebras are fields, an isotropic vector
expresses $c/a$ as a quotient of nonzero norms, and conversely.
If either algebra is split, its norm form is isotropic
and its norm group is $F^\times$, so the equivalence
still holds.

If $r+s\notin\wp(F)$, then either one algebra is split,
or $L_r$ and $L_s$ are distinct quadratic fields.
In the latter case, local class field theory gives
distinct index-two norm subgroups, whose product is
$F^\times$; see \cite{SerreLocalFields}.
Thus $q$ is isotropic in either case.
If $r+s\in\wp(F)$, then $L_r\simeq L_s$, so
$H_rH_s=H_r$, and the cyclic algebra splitting criterion gives
\[
q\text{ is isotropic}
\iff c/a\in H_r
\iff B\text{ splits}.
\]
This proves the assertion.
\end{proof}

We use the splitting certificate $\mathscr S(C,I)$
defined in Section~3.
By Lemma~\ref{lem:local-splitting}, for quaternion algebras
this certificate requires no real splitting hypothesis.

Use the degree-two matched ramification families over $K$,
together with the fixed singleton families at the finite
exceptional places, as in
Proposition~\ref{local-diophantine-correspondence}.
Write their combined parameters as $\lambda\in\Lambda$.
The associated quaternion algebras $R_\lambda$ split at
real places, and the associated uniformly diophantine ideals
satisfy
\[
I_\lambda=\bigcap_{w\in E_\lambda}\Km_w,
\qquad
\varnothing\ne E_\lambda\subseteq\Delta(R_\lambda).
\]
For every finite place $w$, some parameter satisfies
\[
E_\lambda=\{w\},
\qquad I_\lambda=\Km_w.
\]
Here $\Lambda$ is a finite disjoint union of diophantine
parameter spaces, so quantification over $\Lambda$ means
a finite disjunction.
In characteristic two, use the Artin--Schreier families.
In characteristic different from two, use the Kummer families
with $\ell=2$, adjoining the singleton families at the
finite exceptional places.

For such an ideal $I$, define
\[
\Theta(\delta,I)\quad:\Longleftrightarrow\quad
\begin{cases}
\exists z\in K^\times,\ h\in I:
   \delta=z^2(1+4h),&\chara(K)\ne2,\\
\exists z\in K:
   \delta-z^2-z\in I,&\chara(K)=2.
\end{cases}
\]
If $I\subseteq\Km_w$, this implies that $\delta$ is a
nonzero square in $K_w$ in the first case, or belongs to
$\wp(K_w)$ in the second.
Indeed, for $h\in\Km_w$, Hensel's lemma solves
$t^2+t=h$ in $K_w$.
In characteristic different from two,
\[
(1+2t)^2=1+4h,
\]
which also proves the assertion at dyadic places.
In characteristic two, the same equation gives the
Artin--Schreier assertion directly.
If $I=\Km_w$, the converses follow by approximating
a local square root or Artin--Schreier preimage by
an element of $K$.
Thus $\Theta$ is sound at every witness place and
complete when $I=\Km_w$.

\begin{proposition}[A witness for four-dimensional forms]
\label{prop:four-witness}
Let $q,\delta,B$ be as in Lemma~\ref{lem:local-four},
with parameters in a global field $K$.
Then $q$ is anisotropic at some finite completion if and only if
\begin{equation}\label{eq:quadratic-witness}
\begin{gathered}
\exists\lambda\in\Lambda,\ C\in\Cc_2(K):
\quad [C]=[B]+[R_\lambda],\\
\mathscr S(C,I_\lambda),\qquad
\Theta(\delta,I_\lambda).
\end{gathered}
\end{equation}
This condition is diophantine in the parameters of $q$.
\end{proposition}

\begin{proof}
Suppose that the certificate holds, and choose
$w\in E_\lambda$.
Since $I_\lambda\subseteq\Km_w$,
Lemma~\ref{lem:local-splitting} implies that $C_w$ splits,
whereas $(R_\lambda)_w$ is nonsplit.
The equality
\[
[C_w]=[B_w]+[(R_\lambda)_w]
\]
therefore shows that $B_w$ is nonsplit.
The condition $\Theta(\delta,I_\lambda)$ supplies the
remaining hypothesis of Lemma~\ref{lem:local-four},
so $q_w$ is anisotropic.

Conversely, suppose that $q_w$ is anisotropic at a finite
place $w$.
Choose a parameter with
$E_\lambda=\{w\}$ and $I_\lambda=\Km_w$.
By Lemma~\ref{lem:local-four}, $B_w$ is nonsplit.
Hence $[B_w]$ and $[(R_\lambda)_w]$ are the same nonzero
element of $\Br(K_w)[2]$.
By period--index equality, the class $[B]+[R_\lambda]$
has a degree-two representative $C$, taking $C=M_2(K)$
if the class is zero.
Its completion at $w$ splits, so
Lemma~\ref{lem:local-splitting} gives
$\mathscr S(C,\Km_w)$.
The exact-localization property of $\Theta$ gives
$\Theta(\delta,\Km_w)$.

Finally, Brauer-class equality, the splitting certificate,
and $\Theta$ are all diophantine in the displayed parameters.
\end{proof}

\begin{lemma}[Binary parameter tests]
\label{lem:binary-tests}
If $\chara(K)\ne2$ and $d\in K^\times$, then
\[
d\notin K^{\times2}
\quad\Longleftrightarrow\quad
\exists b\in K^\times:
\ (d,b)_{-1}\text{ is nonsplit}.
\]
If $\chara(K)=2$ and $r\in K$, then
\[
r\notin\wp(K)
\quad\Longleftrightarrow\quad
\exists b\in K^\times:
\ [r,b)_2\text{ is nonsplit}.
\]
Both conditions are diophantine over $K$.
\end{lemma}

\begin{proof}
If the corresponding quadratic etale algebra is split,
every associated quaternion symbol splits.
Otherwise, it is a quadratic field extension $L/K$.
Chebotarev provides a finite place $w$ at which $L/K$
has local degree two.
The local norm group is a proper open subgroup of
$K_w^\times$.
By weak approximation, choose $b\in K^\times$ outside
this subgroup.
The corresponding quaternion symbol is then nonsplit
at $w$, hence nonsplit over $K$.
Diophantine definability follows from
Theorem~\ref{prime version}.
\end{proof}

\begin{theorem}\label{thm:general-quadratics}
Let $K$ be a global field and $n\ge1$.
The coefficient tuples
$\mathbf a=(a_{ij})_{1\le i\le j\le n}$ for which
\[
q_{\mathbf a}(X_1,\ldots,X_n)
=\sum_{i\le j}a_{ij}X_iX_j
\]
is anisotropic over $K$ form a diophantine set.
This includes all quadratic forms, without a nonsingularity
assumption, in characteristic two.
\end{theorem}

\begin{proof}
Every change of variables used below is encoded by a matrix
in $\GL_n(K)$ and equality of polynomial coefficients.
It therefore suffices to give diophantine conditions on
the parameters of the indicated normal forms.

Suppose first that $\chara(K)\ne2$.
A degenerate form is isotropic, so it suffices to consider
diagonal forms with nonzero coefficients.
In dimension one, anisotropy is nonvanishing of the coefficient.
In dimension two, $\langle a,b\rangle$ is anisotropic
exactly when $-b/a$ is not a square; apply
Lemma~\ref{lem:binary-tests}.
In dimension three, $\langle a,b,c\rangle$ is anisotropic
exactly when
\[
(-b/a,-c/a)_{-1}
\]
is nonsplit; apply Theorem~\ref{prime version}.

In dimension four, use the first normal form of
Lemma~\ref{lem:local-four}.
By the Hasse--Minkowski principle, anisotropy over $K$
is equivalent to anisotropy at some completion.
The finite-place condition is diophantine by
Proposition~\ref{prop:four-witness}.
At a real place, anisotropy is definiteness, and the
corresponding sign conditions are diophantine by
Lemma~\ref{lem:real-places}.
Complex places give no obstruction.
In dimension at least five, every nondegenerate form is
isotropic at every finite completion, so anisotropy is
equivalent to definiteness at some real place.
For the Hasse--Minkowski principle and the local dimension
bound, see \cite[Theorem~2.3(1)--(2)]{hu2026}.

Now suppose that $\chara(K)=2$.
If $n\ge5$, every quadratic form over $K$ is isotropic
by the $C_2$ property of global function fields;
see \cite[Corollary, p.~378]{Lang1952}.
We may therefore assume that $n\le4$.

We first record that $[K:K^2]=2$.
Indeed, if $\mathbb F_q$ is the constant field and
$[K:\mathbb F_q(t)]=d$, Frobenius gives
$[K^2:\mathbb F_q(t^2)]=d$, whereas
$[K:\mathbb F_q(t^2)]=2d$.
Fix $\tau\in K\setminus K^2$.
Then $\{1,\tau\}$ is a $K^2$-basis of $K$, and
\[
u\notin K^2
\quad\Longleftrightarrow\quad
\exists x,y\in K:
\ u=x^2+\tau y^2,\quad y\ne0.
\]
Thus nonsquareness is diophantine.
Moreover, a totally singular form is a sum of scalar
multiples of squares, so every such form of dimension
at least three is isotropic.

In dimension one, anisotropy means that the coefficient
is nonzero.
In dimension two, a totally singular form
$aX^2+cY^2$ is anisotropic exactly when
$a,c\ne0$ and $c/a\notin K^2$.
A nonsingular binary form has shape $a[1,r]$, with
$a\ne0$, and is anisotropic exactly when $r\notin\wp(K)$.
Both conditions are diophantine, the latter by
Lemma~\ref{lem:binary-tests}.

In dimension three, a form of polar rank zero is totally
singular and hence isotropic.
A form of polar rank two has shape
\[
a[1,r]\perp cZ^2,
\qquad a\ne0.
\]
If $c=0$, it is isotropic.
If $c\ne0$, it is anisotropic exactly when
$[r,c/a)_2$ is nonsplit, by the norm criterion.
Thus Theorem~\ref{prime version} applies.

In dimension four, first suppose that the polar form
is singular.
Polar rank zero has already been treated.
In polar rank two, decompose the form into a nonsingular
binary form and its two-dimensional polar radical.
If the restriction to the radical is isotropic, so is
the whole form.
Otherwise, its two coefficients form a $K^2$-basis of $K$,
so the radical represents every element of $K$.
It can therefore cancel the value of a nonzero vector
in the other summand.
Thus every four-dimensional form with singular polar
form is isotropic.

It remains to consider nonsingular four-dimensional forms.
Such a form has shape
\[
a[1,r]\perp c[1,s],
\qquad a,c\in K^\times,\quad r,s\in K.
\]
Indeed, decompose its polar form into symplectic planes,
choose a vector of nonzero quadratic value in each plane,
and normalize its partner.
By the Hasse--Minkowski principle
\cite[Theorem~3.2]{Pollak1970}, anisotropy over $K$
is equivalent to anisotropy at some completion.
All places are nonarchimedean, so
Proposition~\ref{prop:four-witness} gives the required
diophantine condition.

In each characteristic, the anisotropic forms are
therefore described by finitely many normal-form families
with diophantine parameter conditions.
Taking their finite union and projecting the corresponding
change-of-variables relations proves the assertion.
\end{proof}
\begin{remark}
In characteristic different from two, the conclusion is
\cite[Theorem~1.3]{hu2026}; the argument above gives an
alternative proof using the present localization families.
In characteristic two, $\Theta$ ensures that the Arf
invariant vanishes at the finite witness place, without
requiring it to vanish over $K$.
Singular forms are treated separately using $[K:K^2]=2$,
without invoking a local--global principle for them.
\end{remark}
\subsection{Generalized Severi--Brauer varieties}

For a central simple algebra $A$ of degree $n$, let
$\SB_r(A)$ denote the generalized Severi--Brauer variety
parametrizing right ideals of reduced dimension $r$,
that is, of $K$-dimension $nr$.

\begin{proposition}\label{prop:generalized-SB}
Let $K$ be a global field and $1\le r<n$.
The set
\[
\{\mathbf c\in\Cc_n(K):
  \SB_r(A_{\mathbf c})(K)=\varnothing\}
\]
is diophantine over $K$.
\end{proposition}

\begin{proof}
Write $A\simeq M_m(D)$, where $D$ is a central division
algebra of degree $d=\ind(A)$, so $n=md$.
Under the identification $A=\End_D(D^m)$, right ideals
correspond to right $D$-subspaces $W\subseteq D^m$ via
\[
W\longmapsto
\{f\in\End_D(D^m):\operatorname{im}(f)\subseteq W\}.
\]
If $\dim_D W=j$, the corresponding ideal has
$K$-dimension $mjd^2=njd$, hence reduced dimension $jd$.
Consequently,
\[
\SB_r(A)(K)\ne\varnothing
\quad\Longleftrightarrow\quad
\ind(A)\mid r.
\]

Over a global field, period equals index.
Therefore
\begin{equation}\label{eq:SB-tensor}
\SB_r(A)(K)=\varnothing
\quad\Longleftrightarrow\quad
r[A]\ne0
\quad\Longleftrightarrow\quad
A^{\otimes r}\text{ is nonsplit}.
\end{equation}
The algebra $A^{\otimes r}$ has fixed degree $n^r$,
and its structure constants are polynomial in those of $A$.
The conclusion follows from Theorem~\ref{main theorem}.
\end{proof}

\begin{proof}[Proof of Corollary~\ref{cor:specialization}]
Choose a $K(T)$-basis of $A$ beginning with $1$, and let
$c_{ij}^k(T)\in K(T)$ be its structure constants. Let $R(T)$ be the
matrix of the canonical central-simplicity map used in
Lemma~\ref{csa locus}. Since $A$ is central simple,
$\det R(T)\ne0$. Choose a nonzero polynomial $h\in K[T]$ divisible
by the denominators of the structure constants and by a numerator of
$\det R(T)$. Put $U=\{t:h(t)\ne0\}$.

The tables over $K[T,h^{-1}]$ define a finite free algebra for which
the canonical central-simplicity map is invertible; equivalently, an
Azumaya algebra. Every specialization at $t\in U(K)$ is therefore
central simple of degree $n$. The graph of
$t\mapsto(c_{ij}^k(t))$ on $U(K)$ is diophantine: impose
$h(t)z=1$ and clear the fixed denominators in the equations defining
the graph. Pull back the diophantine set of
Proposition~\ref{prop:generalized-SB} along this graph.
\end{proof}
\section*{Funding}
This work was carried out during a research visit to TU Dresden
supported by the Sino-German (CSC--DAAD) Postdoc Scholarship Program.
The work was also supported by the Shandong Provincial Natural Science
Foundation (Project No.~ZR2024QA184).
\section*{AI Disclosure}
During the preparation of this manuscript, the author used OpenAI's
GPT-5.6 Sol as a research aid. In particular, although the author had
long suspected the existence of the mechanism underlying the local
witness argument, the model proposed a proof strategy that helped
formulate it precisely. AI assistance was also used in drafting
auxiliary proofs and in developing the revised organization and the
Kummer argument. The author is responsible for verifying all statements,
proofs, and references, and for the entire content of the manuscript.

\bibliographystyle{alpha}
\bibliography{SBV}
\end{document}